\documentclass{amsart}

\usepackage{amsfonts}
\usepackage{latexsym}
\usepackage{amsmath}
\usepackage{epsfig}
\usepackage{amssymb}
\usepackage{enumerate}
\usepackage[colorlinks=true,citecolor=blue,linkcolor=blue,urlcolor=blue,backref=true,pagebackref=true,pdfauthor={Philippe Gambette, Katharina T. Huber},pdftitle={A Note on Encodings of Phylogenetic Networks of Bounded Level}]{hyperref}

\usepackage[leqno]{amsmath}
\usepackage{amscd}
\usepackage{epsfig}

\newtheorem{lemma}{Lemma}
\newtheorem{corollary}{Corollary}

\newtheorem{theorem}{Theorem}
\newtheorem{proposition}{Proposition}

\newcommand{\old}[1]{}

\newcommand{\cR}{{\mathcal R}}

\newcommand{\pf}{\noindent{\em Proof: }}
\newcommand{\epf}{\hfill\hbox{\rule{3pt}{6pt}}\\}

\begin{document}

\title{A Note on Encodings of Phylogenetic Networks of Bounded Level}
\thanks{
This work was supported by the French ANR projects ANR-06-BLAN-0148-01
(GRAAL) and ANR-08-EMER-011-01 (PhylARIANE), and was initiated during
the MIEP workshop in 2008.
%~\\
%Philippe Gambette, D{\'e}partement Informatique, L.I.R.M.M, Universit{\'e} Montpellier II, C.N.R.S., France.\\
%Katharina T. Huber, School of Computing Sciences, University of East Anglia, Norwich, NR4 7TJ, UK.
}

%\titlerunning{A Note on Encodings of Level-k Networks}

\author{Philippe Gambette, Katharina T. Huber}

\date{\today}

\maketitle

\begin{abstract}
Driven by the need for better models that allow one to
 shed light into the question
how life's diversity has evolved, phylogenetic networks 
%(essentially leaf-labelled rooted directed acyclic graphs) 
have now joined phylogenetic trees 
%(essentially leaf labelled rooted graph-theoretical trees) 
in the center of phylogenetics research. Like phylogenetic trees, 
such networks canonically induce collections of phylogenetic trees, 
clusters, and triplets, respectively. Thus it is not surprising that
many network approaches aim to reconstruct a phylogenetic network from
such collections. Related to the well-studied perfect
phylogeny problem, the following question is of fundamental importance in this context: 
When does one of the above collections
encode (i.e. uniquely describe) the network that induces it?

In this note, we present a complete answer to this question 
for the special case of a level-$1$ (phylogenetic) network by 
characterizing those level-$1$ networks
for which an encoding in terms of one (or equivalently all)
of the above collections exists. Given that this type of network
forms the first layer of the rich hierarchy of level-$k$ networks,
$k$ a non-negative integer, it is natural to wonder whether 
our arguments could be extended
to members of that hierarchy for higher values for $k$. By giving
examples, we show that this is not the case.

%Include keywords, PACS and mathematical
%subject classification numbers as needed.
\end{abstract}

\section{Introduction}
\label{intro}

%\subsection{Phylogenetic Networks}

An improved understanding of the complex processes that drive evolution
has lent support to the idea that reticulate evolutionary events such as
lateral gene transfer or hybridization are more common than originally
thought rendering a phylogenetic tree
(essentially a rooted leaf labelled graph-theoretical tree) 
too simplistic a model to fully understand the complex processes 
that drive evolution. 
Reflecting this, phylogenetic networks have now
joined phylogenetic trees in the center of phylogenetics research.
Influenced by the diversity of questions posed by evolutionary 
biologists that can be addressed with a phylogenetic networks, various
alternative definitions of these types of networks have been developed
over the years~\cite{HusonBryant2006}. These include split networks
\cite{BryantMoulton2004,BFSR1995,HHML2004} 
as well as ancestral recombination graphs \cite{SongHein2005}, TOM 
networks \cite{Willson2006a},
level-$k$ networks\footnote{Note that 
these networks were originally introduced in \cite{JanssonSung2004b},
but the definition commonly used now is slightly different with the 
main difference being that every vertex of the network with indegree
2 must have outdegree 1 (see e.g.
\cite{IKKS2008} and the references therein).}
%\cite{JanssonSung2004b}
with $k$ a non-negative integer that
in a some sense captures how complex the network structure is,
networks for studying the evolution of polyploid organisms
\cite{HuberMoulton2006}, tree-child and tree-sibling 
networks \cite{CLRV2008b}, to name just a few. 

\begin{figure}[!hgt]
\centering
\begin{tabular}{ccc}
(a) & (b) & (c) \\
\includegraphics[scale=0.45]{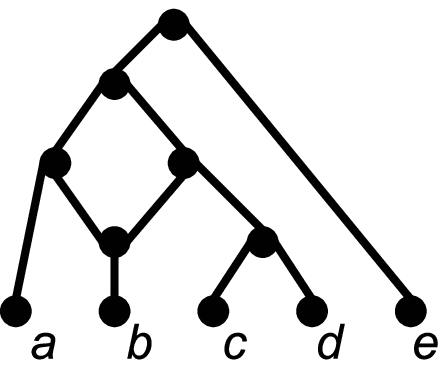} &
\includegraphics[scale=0.45]{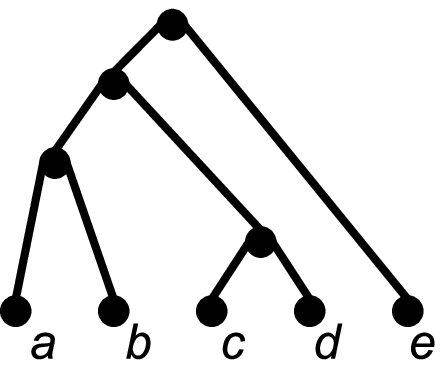} &
\includegraphics[scale=0.45]{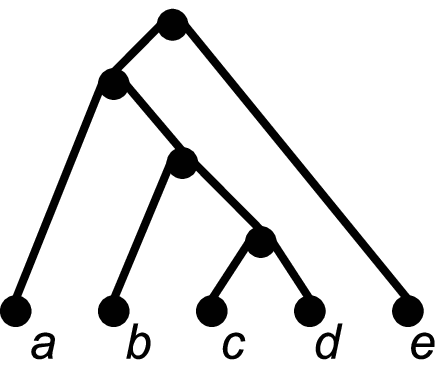}\\
%\includegraphics[scale=0.45]{Level1NetworkNonMinimal.eps} &
%\includegraphics[scale=0.45]{Level1NetworkOther.eps} &
%\includegraphics[scale=0.45]{Level1NetworkMultree.eps} \\
%(d) & (e) & (f)
\end{tabular}
\caption{(a) A level-$1$ phylogenetic network $N$. (b) 
and (c) The  phylogenetic trees that form the tree system
$\mathcal T(N)$.
}
\label{FigLevel1}
\end{figure}

Apart from split networks which aim 
to give an implicit model of evolution and are not the focus of this
note, all other phylogenetic networks mentioned above aim to provide an 
explicit model of evolution. Although slightly different in detail, they 
are all based on the concept of a leaf-labelled 
rooted connected directed acyclic graph (see 
the next section for a definition).
For the convenience of the reader, we depict an 
example of a phylogenetic network in the form of a
level-$1$ network in Fig.~\ref{FigLevel1}(a).  
Concerning these types of phylogenetic networks, it should be noted that 
they are closely related to 
{\em galled trees}~\cite{WZZ2001,GEL2003} and that, in addition to constituting 
the first layer of the rich hierarchy of level-$k$ networks, they also
form a subclass of the large class of 
tree-sibling networks \cite{AVP2008}.

Due to the rich combinatorial structure of phylogenetic networks,
different combinatorial objects have been used to reconstruct them from
biological data. For a set $X$ of taxa (e.g. species or organisms), these
include {\em cluster systems} of $X$, 
that is, collections of subsets
of $X$ \cite{BandeltDress1989,HusonRupp2008},
 {\em triplet systems} on $X$, that
is, collections of phylogenetic trees with just three leaves
which are generally called {\em (rooted) triplets}~\cite{JanssonSung2004b,ToHabib2009},
and {\em tree systems}, that is,
collections of phylogenetic trees which all have leaf 
set $X$~\cite{Semple2007}.
The underlying rational being that any phylogenetic network
$N$ induces a cluster system $\mathcal C(N)$, a triplet
system $\mathcal R(N)$ and a tree system
$\mathcal T(N)$. Again we defer the precise definitions to 
later sections of this note and  remark that for the
level-$1$ network $N$ with leaf set $X=\{a,b\ldots,e\}$
depicted in Fig.~\ref{FigLevel1}(a),
the cluster system $\mathcal C(N) $ consists of $X$,
the five singleton sets of $X$, and the subsets 
$\{a,b\}, \{c,d\}, \{b,c,d\}, Y:=\{a,b,c,d\}$, and
the tree system $\mathcal T(N)$ consists of the
phylogentic trees depicted in
Fig.~\ref{FigLevel1}(b) and (c), respectively. Denoting a
phylogenetic tree $t$ on $x,y,z$ such that the root of 
$t$ is not the parent vertex of $x$ and $y$ by $z|xy$ (or equivalently
by $xy|z$) 
then the triplet system $\mathcal R(N)$ consists
of all triplets of the from $e|xy$ where $x,y\in Y$ distinct, 
$x|cd$ with $x\in \{a,b\}$, and
$x|ab$ and $a|bx$ with
$x\in \{c,d\}$.
%
%The \emph{softwired cluster system} $\mathcal S(N)$ is the union of
%the cluster systems of all the trees of $\mathcal T(N)$, which happens
%to be equal to $\mathcal C(N)$ for this example.

Although undoubtedly highly relevant for 
phylogenetic network reconstruction, 
the following fundamental question has however remained 
largely unanswered so far (the main exception being the case
when $N$ is in fact a phylogenetic tree in which
case this question is closely related to the well-studied 
{\em perfect phylogeny problem} -- see e.g.~\cite{GruenewaldHuber2007} for a recent overview.): 
When do the systems 
$\mathcal C(N)$, $\mathcal R(N)$, or $\mathcal T(N)$ induced by
a phylogenetic network $N$ {\em encode} $N$, that is, there is no
other phylogenetic network $N'$ for which the corresponding
systems for $N$ and $N'$ coincide? 

Complementing the insights for when $N$ is a phylogenetic tree alluded to above,
 answers were recently provided for $\mathcal R(N)$ 
in case $N$ is a very special type of level-$k$ network, $k\geq 2$, 
\cite{IKM2009} and for $\mathcal T(N)$ for the special case that $N$ 
is a regular network \cite{Willson2009}.
Undoubtedly important first results, there are many types of phylogenetic 
networks which are encoded by the tree system
they induce but which are not regular or by the triplet system
they induce but do not belong to  that special  class of
level-$2$ networks. An example for both cases is the level-$1$ network 
depicted in Fig.~\ref{FigLevel1}(a). Although one might be tempted to
speculate that all level-$1$ networks enjoy this property, this is not the
case since the level-$1$ networks depicted in Fig.~\ref{FigLevel1}(a) and 
Fig.~\ref{Figminimal}(b),
respectively, induce the same tree system and the same triplet system. 
%(and also the same triplet systems). 
The main result of this paper shows that 
these observations are not a coincidence.
More precisely, in Theorem~\ref{theorem-montpellier} we 
 establish that a level-$1$ network $N$
is encoded by the triplet system $\mathcal R(N)$ (or equivalently
by the tree system $\mathcal T(N)$ or equivalently
the  cluster system $\mathcal S(N)=\mathcal S({\mathcal T}(N))
:=\bigcup_{T\in \mathcal T(N)}\mathcal C(T)$
which arises in the context of the {\em softwired interpretation} 
of $N$ \cite{HusonRupp2008} and contains $\mathcal C(N)$) 
if and only if, when ignoring directions, $N$ does not
contain a cycle of length 4. Consequently the number 
of non-isomorphic (see below) phylogenetic networks $N'$ which all
induce the same tree system (or equivalently
 the same triplet system or the same cluster system $\mathcal S(N)$) 
 grows exponentially in the number of cycles of $N$ of length 4.
It is of course highly tempting to speculate that a similar
characterization might hold for higher values of $k$. 
However as our examples
show, establishing such a result will require an alternative approach since
 our arguments cannot be extended to level-$2$ networks
and thus to level-$k$ networks with $k\geq 2$.

This note is organized as follows. In the next section, we present 
the definition of a level-$1$ network plus surrounding terminology. 
In Section~\ref{clusterEncoding}, we
present the definitions of the cluster system $\mathcal C(N)$ and
the tree system $\mathcal T(N)$ induced by a phylogenetic network $N$.
This also completes the definition of the cluster system $ \mathcal S(N)$
given in the introduction.
Subsequent to this, we show that for any
level-$1$ network $N$, the cluster systems $ \mathcal S(N)$
and $\mathcal C(N)$ are
{\em weak hierarchies} (Proposition~\ref{softwired-weak-hierarchy})
which are well-known in cluster analysis. 
In addition, we show that this property is not enjoyed by
level-$2$ networks and thus level-$k$ networks with $k\geq 2$.
%Further, we present a necessary condition
%for a cluster system to be induced by a level-$1$ network
%when interpreted as a softwired level-$1$ network.
In Section~\ref{Simple-level-$1$}, we first present the definition of the
triplet system $\mathcal R(N)$ induced by phylogenetic network
$N$. Subsequent to this, we turn our attention to 
the special case of encodings of 
simple level-$1$ networks. In Section~\ref{alternativeEncoding}, 
we present our main result (Theorem~\ref{theorem-montpellier}). 

To ease the presentation of our results, 
%from now on all rooted pylogenetic trees and
%networks will be called pylogenetic trees and networks. Also
%from now on $X$ denotes a finite set. Also in a
in all figures the (unique) root of a
network is the top vertex and all arcs are directed downwards, 
away from the root. Furthermore, for any directed graph $G$, we denote the
vertex set of $G$ by $V(G)$ and the set of arcs of $G$ by  $A(G)$.

%\subsection{Outline of the Article}

\section{Basic terminology and results concerning level-$1$ networks} 
\label{prelim}

In this section we present the definitions
of a phylogenetic network and of a level-$k$ network, $k\geq 0$.
In addition we also provide  the basic terminology 
surrounding these structures. 

Suppose $X$ is a finite set.
A {\em  phylogenetic network} $N$ on $X$ is a rooted 
directed acyclic graph (DAG) that satisfies the following additional
properties. (i) The set $L(N)$ of {\em leaves} 
of $N$, that is vertices with indegee
1 and outdegree 0, is $X$.
(ii) Exactly one vertex of $N$, called the
{\em root}  and denoted by $\rho_N$, has indegree 0 and outdegree 2.
(iii) All vertices of $N$ that are not contained in $L(N)\cup \{\rho_N\}$
are either {\em split vertices}, that is, have indegree 1 and outdegree 2 
or {\em reticulation vertices}, that is, have indegree 2 and outdegree 1. 
The set of reticulation vertices of $N$ is denoted by $R(N)$. 
A phylogenetic network $N$ for which
$R(N)$ is empty is called a {\em (rooted) phylogenetic tree (on $X$)}.
Two phylogenetic networks $N$ and $N'$ which both have
leaf set $X$ are said to be {\em isomorphic} if there exists a bijection
from $V(N)$  to $V(N')$  which is the identity on $X$ and 
induces a graph isomorphism between $N$ and $N'$.

%A {\em (rooted) triplet} $t$ with vertex set $Y=\{x,y,z\}$ 
%is a phylogenetic tree on $Y$. If the parent of
%$\{y,z\}$ in $t$ is not the root of $t$ then we 
% denote $t$  by $x|yz$ (or alternatively by $x|zy$ or $yz|x$). 

To present the definition of a level-$k$ network, we need to
introduce some terminology concerning rooted DAGs first.
Suppose $G$ is a rooted connected DAG with at least 2 vertices. 
%Then we call an arc 
%$e\in A(G)$ a {\em cut-arc} of $G$ if the deletion of $e$ 
%from $G$ disconnects $G$ and a {\em non-cut-arc} 
%otherwise. Also for any path $P$ in $G$ from a vertex $a\in V(G)$
%to a vertex $b\in V(G)$ we call the vertices in $V(P)-\{a,b\}$
%the {\em interior} vertices of $P$. 
Then we denote the graph obtained from  
$G$ by ignoring the directions on $G$ by $U(G)$. If $H$ is a graph with
at least 2 vertices  then we call
$H$  \emph{biconnected} if $H$ does not contain a vertex
whose removal disconnects it. A {\em biconnected component} of $H$
is a maximal subgraph of $H$ that is biconnected. If $G$ is a phylogenetic
network and $B$ is a rooted sub-DAG such that $U(B)$
is a biconnected component of $U(G)$ then we call $B$ a {\em blob}.

\begin{figure}[!hgt]
\centering
\begin{tabular}{ccc}
(a) & (b)  \\
\includegraphics[scale=0.45]{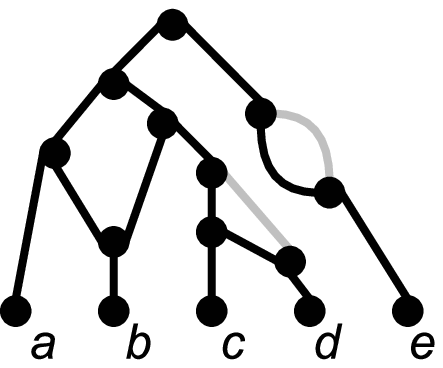} &
\includegraphics[scale=0.45]{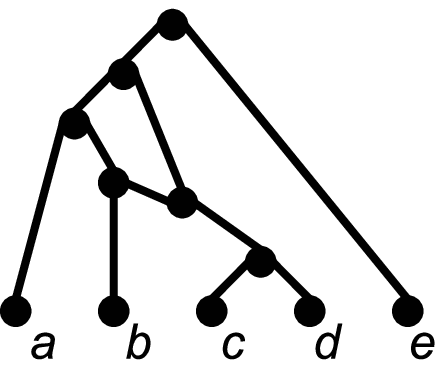} &
\end{tabular}
\caption{
The level-$1$ network $N$ depicted in (a) induces and thus
represents the same triplet system $\mathcal R(N)$, cluster systems
$\mathcal C(N)$ and $\mathcal S(N)$, and tree system $\mathcal T(N)$
as the level-$1$ network $N'$ presented in (b). However,  $N'$ is 
a less parsimonious representation of these 4 systems. 
}
\label{Figminimal}
\end{figure}

Following~\cite{IKKS2008}, we call a phylogenetic network $N$ a
{\em level-k} network for some non-negative integer $k$
if each blob of $N$ contains at most
$k$ reticulation vertices. Note that some authors 
define a level-$1$ network $N$ to be a phylogenetic network
without the above outdegree requirement on the reticulation vertices
of $N$ (see e.g.~\cite{JanssonSung2004b}). Also and sometimes on its own
or in addition to the above, the requirement
that each blob contains at most $k$ reticulation vertices is 
sometimes replaced by the requirement that the cycles in $U(N)$ are node
disjoint (see e.g.~\cite{JanssonSung2004b,JanssonSung2006}). 
Although in spirit the same definitions, the difference is that a cycle is generally 
understood to have at least three vertices which implies that the
network depicted in Fig~\ref{Figminimal}(a) would not be a 
level-$1$ network. However 
the definition presented in~\cite{IKKS2008} would render
that network a level-$1$ network. Having said that, 
the network $N$ depicted in Fig.~\ref{Figminimal}(b) is 
a less parsimonious representation of the same biological information
(expressed in  terms of the systems $\mathcal T(N)$,
$\mathcal R(N)$, $\mathcal C(N)$, and $\mathcal S(N)$)
as the level-$1$ network in  Fig.~\ref{Figminimal}(a) in the sense that
the edges in grey are redundant for displaying that information. To avoid
these types of level-$1$ networks which cannot be encoded by any of the
4 systems of interest in this note, we follow~\cite{IKM2009}
and require that every blob in a level-$1$ network
contains at least 4 vertices.
%
%In addition, we call a level-$k$
%network $N$, {\em strict} if either $k=0$ (i.e.
 %$N$ is a phylogenetic tree) or $k\geq 1$ and $N$ is 
%not also a level-$(k-1)$ network.
%Note that although deceptively simple looking
%for $k=1,2$, these networks can give rise to surprisingly
%difficult algorithmic problems (see 
%e.\,g.\,\cite{JanssonSung2006,IKKS2008,ToHabib2009}.

\begin{figure}[!hgt]
\centering
\begin{tabular}{cccccc}
\includegraphics[scale=0.5]{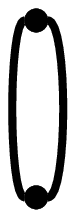}
& ~ ~ ~
& \includegraphics[scale=0.5]{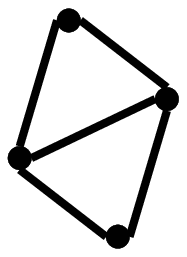}
& \includegraphics[scale=0.5]{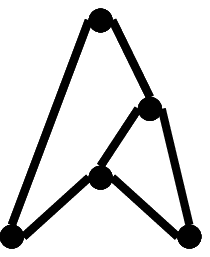}
& \includegraphics[scale=0.5]{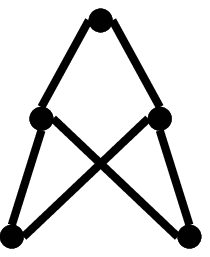}
& \includegraphics[scale=0.5]{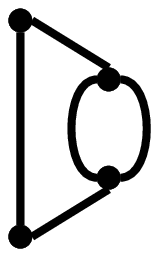}\\
$\mathcal{G}^1$ & & $\mathcal{G}^2_a$ & $\mathcal{G}^2_b$ & 
$\mathcal{G}^2_c$ & $\mathcal{G}^2_d$
\end{tabular}
\caption[level-$1$ and level-$2$ generators]{The unique 
level-$1$ generator $\mathcal{G}^1$, 
and the four level-2 generators:
$\mathcal{G}^2_a$, $\mathcal{G}^2_b$, $\mathcal{G}^2_c$ and $\mathcal{G}^2_d$.}
\label{level-k-generators}
\end{figure}

For $k=1,2$, it was shown in~\cite{IKKS2008} (see also~\cite{JanssonSung2006} for the
case $k=1$) that level-$k$ networks can be built up by chaining together
structurally very simple level-$k$ networks called {\em simple level-$k$ networks}.
Defined for general non-negative integers $k$, these atomic building blocks
are precisely those level-$k$ networks that can 
be obtained from a {\em level-$k$ generator} by
applying a certain ``leaf hanging'' operation \cite{IKKS2008} to its
``sides''. Such a generator is a biconnected directed acyclic multi-graph 
which has a single root, precisely $k$ {\em pseudo-reticulation
vertices} (i.\,e.\,vertices with indegree 2 and outdegree
at most 1) and all other vertices are split vertices where the
root and a split vertex are defined as in the case 
of a phylogenetic network. For the convenience of the reader,
we present in Fig~\ref{level-k-generators} the
unique level-$1$ generator and all 4 level-$2$ generators 
which originally appeared in slightly different form
in~\cite{IKKS2008}.
Regarding larger values for $k$, it was recently shown 
in~\cite{L3Generators}
that there exist $65$ level-3 generators. In addition, it
was shown in~\cite{GBP2009}
that there are $1993$ level-4 generators and that the number
of level-$k$ generators grows exponentially in $k$.
A {\em side} of a generator $G$ is an arc of $G$ or one of its 
pseudo-reticulation vertices.

From now on and unless stated otherwise, all phylogenetic 
networks have leaf set $X$.

\section{The Systems $\mathcal C(N)$, $\mathcal T(N)$, and 
$\mathcal S(N)$ }\label{clusterEncoding}

In this section, we introduce for a phylogenetic
network $N$ the associated systems $\mathcal C(N)$, 
$\mathcal T(N)$, and $\mathcal S(N)$ already mentioned in the
introduction. In addition, we prove that 
in case $N$ is a level-$1$ network the associated
systems $\mathcal C(N)$ and $\mathcal S(N)$ are weak hierarchies. 
We conclude with presenting an example
that shows that level-$k$ networks, $k\geq 2$, do not enjoy this 
property in general. We start with some definitions.

Suppose $N$ is a phylogenetic network.
Then we say that a vertex $a\in V(N)$ is {\em below}
a vertex $b\in V(N)$  denoted by $a \preceq_N b$, if there 
exists a path $P_{ba}$
(possibly of length 0) from $b$ to $a$.  In this case, we also say
that $b$ is {\em above $a$}. Every vertex $v\in V(N)$
therefore induces a non-empty subset $C(v) =C_N(v)$ of $ X$  
which comprises of all leaves of $N$ below $v$ (see e.g. 
\cite{SempleSteel2003}). 
We collect the subsets $C(v)$ induced by the vertices $v$ of $N$
this way in the set $\mathcal C(N)$, i.e. we put 
$\mathcal C(N)=\bigcup_{v\in V(N)}\{C(v)\} $. 
For convenience, we refer to any collection $\mathcal C$ of 
non-empty subsets of
$X$ as a {\em cluster system (on $X$)} and to the elements of $\mathcal C$
as {\em clusters} of $X$.
It should be noted that in case $N$ is a binary phylogenetic tree,
the cluster system $\mathcal C(N)$ is a {\em hierarchy (on $X$)},
that is, for any two clusters $C_1,C_2\in \mathcal C(N) $
we have that $C_1\cap C_2\in \{\emptyset, C_1,C_2\}$. Hierarchies
are sometimes also called {\em laminar families}, and it is well-known
that the set of clusters $\mathcal C(T)$
induced by a binary phylogenetic tree $T$  
uniquely determines that tree (see e.g. \cite{SempleSteel2003}).

In the context of phylogenetic network construction, the concept 
of a {\em weak hierarchy (on $X$)} was introduced
in \cite{BandeltDress1989}. These objects are defined as follows.
Suppose $\mathcal C$ is a cluster system on $X$. Then  
$\mathcal C$ is called a {\em weak hierarchy (on $X$)} if 
\begin{eqnarray}\label{helly}
C_1\cap C_2\cap C_3\in \{C_1\cap C_2, C_2\cap C_3, C_1\cap C_3\}
\end{eqnarray}
holds for any three
elements $C_1,C_2, C_3\in {\mathcal C}$.
Note that the above property is sometimes 
also called the {\em weak Helly property} \cite{SempleSteel2003}. 
Also note that any hierarchy is in particular a weak hierarchy and that
any subset of a weak hierarchy is again a weak hierarchy. Finally note
that weak hierarchies are
well-known objects in classical hypergraph and abstract
convexity theories \cite{BandeltDress1989} 
(see also the reference therein and \cite{BBO04}), and that they
where originally introduced into cluster analysis as {\em medinclus}
in \cite{B88,B89}. 

We will establish the main result of this 
section (Proposition~\ref{softwired-weak-hierarchy}) by showing
that the cluster system $\mathcal S(N)$ associated to a
level-$1$ network $N$ is a weak hierarchy. 
To do this, we first
need to complete the definition of $\mathcal S(N)$
which relies on the definition of the system $\mathcal T(N)$. 
We will do this next. 

Suppose $N$ is a phylogenetic network. Then  we say that
a phylogenetic tree $T$ is {\em displayed} by $N$ if 
the leaf set of $T$ is $X$ and $T$ is a phylogenetic
tree  obtained from $N$ via the following process. 
For each reticulation vertex of $N$ delete 
one incoming arc and suppress any resulting 
degree 2 vertices. In case the root $\rho_N$ of $N$ is rendered 
a vertex with out-degree 1 this way, we  identify $\rho_N$ with its unique 
child. The set $\mathcal T(N)$ is the collection of all 
phylogenetic trees that are displayed by $N$. 
%Following \cite{HusonRupp2008},
%we will sometimes refer to a phylogenetic network $N$  as a
%{\em softwired phylogenetic network} when we view it in terms of  the
%tree system $\mathcal T(N)$ that $N$ induces.
%%Combined with the definition outline for $\mathcal S(N)$ given in the
%%introduction, it follows that
%%$\mathcal S(N)$ is the collection of all clusters $C$ on $L(N)$
%%for which there exists a phylogenetic tree $T\in \mathcal T(N)$ such
%%that $C\in\mathcal C(T)$. 
To every vertex $v\in V(N)$
a cluster system  $\mathcal S_N(v)$ defined by putting 
\begin{align*}
\mathcal S_N(v) =  \{C_T(v): T \in \mathcal T(N)\}
\end{align*}
can be associated. Clearly, $C_N(v) \in \mathcal S_N(v) $ and
$\mathcal S(N)=\bigcup_{v\in V(N)} \mathcal S_N(v)$. 

To link clusters of $X$ with level-$1$ networks on $X$,
we say that a cluster $C$ on $X$  is {\em tree-consistent with a 
level-$1$ network} $N$ if $C \in \mathcal S(N)$.
More generally, we say that a cluster system $\mathcal C$ is
{\em tree-consistent with a level-$1$ network}
$N$ if $\mathcal C\subseteq \mathcal S(N)$ holds.
Thus, for any level-$1$ network $N$ the cluster system $\mathcal S(N)$ 
equals the set of all clusters of $X$ that are
tree-consistent with $N$.
Finally, we say that a cluster
system $\mathcal C$ is {\em level-$1$-consistent} if there
exists a level-$1$ network
$N$ such that $\mathcal C$ is tree-consistent with $N$.

We next establish Proposition~\ref{softwired-weak-hierarchy}.
Its proof relies on a characterization of a weak hierarchy $\mathcal H$
on $X$  from \cite[Lemma 1]{BandeltDress1989} in terms of a property of
a certain  $\mathcal H$-closure that can be canonically associated to 
$\mathcal H$. More precisely, suppose  $\emptyset\not=Y\subseteq X$ 
and $\mathcal H$ is a cluster system on $X$. Then the 
{\em $\mathcal H$-closure} 
$\langle Y\rangle_{\mathcal H}$
of $Y$ is the intersection $\bigcap_{Y\subseteq C,\,\, C\in \mathcal H}C$.
%where we put $X=\bigcap_{\emptyset\subseteq C,\,\, C\in \mathcal H}
%\emptyset$. 
Now a cluster system $\mathcal H$ on $X$ is a weak hierarchy if and  
only if for every non-empty subset $A\subseteq X$  
there exists elements $a,a'\in A$ such that 
$\langle A\rangle_{\mathcal H}=\langle \{a,a'\}\rangle_{\mathcal H}$. 
Note that this implies in particular
that the number of elements in a weak hierarchy is at most ${|X|+1}\choose{2}$
\cite{BandeltDress1989}. With regards to this bound it should be
noted that it was recently shown in \cite{KNTX2008}
that the size of a cluster system which is tree-consistent 
with a level-$1$ network $N$ is linear in $|X|$.
In view of Proposition~\ref{softwired-weak-hierarchy}, this bound 
improves on the previous  bound for this special kind of weak 
hierarchy.

\begin{proposition}\label{softwired-weak-hierarchy}
A level-$1$-consistent cluster system
is a weak hierarchy. In particular, the systems  $\mathcal S(N)$
and $\mathcal C(N)$ associated to a level-$1$ network $N$ are 
weak hierarchies. 
\end{proposition}
\pf Since 
every subset of a weak hierarchy is again a weak hierarchy, 
it suffices to show that for every level-$1$ network $N$ the 
associated cluster system  $\mathcal S(N)$ is a weak hierarchy.
To see this suppose $N$ is a level-$1$ network on $X=\{x_1, \ldots, x_n\}$,
$n\geq 1$. 
Consider a graphical representation of $N$ and, starting from the left most
leaf of $N$ in that representation, let $x_1\ldots x_n$ denote 
the induced ordering of the leaves of $N$ (note that this might involve
re-labelling some of the elements in $X$). Suppose  
$\emptyset \not=A\subseteq X$. Let $i,j\in\{1,\ldots, n\}$ be such that 
$x_j\in A$  and every  leaf in
$X$ succeeding $x_j$ in that ordering is not contained in $A$. Similarly,
let $x_i\in A$ be such that every leaf in
$X$ preceding $x_i$ in that ordering is not contained in $A$. We claim
that 
$\langle A\rangle_{\mathcal S(N)}=\langle \{x_i,x_j\}\rangle_{\mathcal S(N)}$.
To see this, note that since $N$ is a level-$1$ network, there exists 
a subtree $T$ of $N$ such that the leaf set of $T$ is $A$. Note that
$T$ might contain vertices whose indegree and outdegree is one. By 
deleting for each reticulation vertex below the root of
$T$ one of its incommming arcs and supressing the resulting degree
2 vertex $T$ can be canonically extended to a subtree $T'$ of some 
tree $T''\in\mathcal T(N)$
such that $\{x_i,x_{i+1},\ldots, x_j\} \subseteq L(T')$ and  
$L(T')$ is minimal with regards to set inclusion. Note that  
$L(T')\in \mathcal S(N)$. But then, by construction,
$\langle A\rangle_{\mathcal S(N)}= L(T')=
\langle \{x_i,x_j\}\rangle_{\mathcal S(N)}$ which proves the claim.
\epf

We remark in passing that to any cluster system $\mathcal C$ of $X$
a similarity 
measure $D_{\mathcal C}:X\times X\to \mathbb R$ 
can be associated to $\mathcal C$ by putting 
$D_{\mathcal C}(a,b)=|\{C\in \mathcal C: a,b\in C\}|$, $a,b\in X$. 
Proposition~\ref{softwired-weak-hierarchy} 
combined with the main result from \cite{BandeltDress1989}
implies that any tree-consistent cluster system $\mathcal C$
can be uniquely reconstructed from its associated similarity measure
$D_{\mathcal C}$. Using the well-known {\em Farris transform} 
(see e.\,g.\,\cite{SempleSteel2003}, and \cite{DHM2007} for a recent 
overview) a similarity
measure can be canonically transformed into a distance measure
$D^{\mathcal C}$ on $X$, that is, a map on $X\times X$ into the non-negative reals 
that is symmetric, satisfies the triangle inequality, and vanishes on the
main diagonal. The latter measures were recently investigated in 
\cite{CJLY2005} from an algorithmical point of view in the context of 
representing them in terms of an {\em ultrametric} level-$1$
network. These are generalizations of ultrametric phylogenetic trees
in the sense that every path from the root of the network 
to any leaf is of the same length.

\begin{figure}[!hgt]
\centering
\begin{tabular}{cccc}
\includegraphics[scale=0.3]{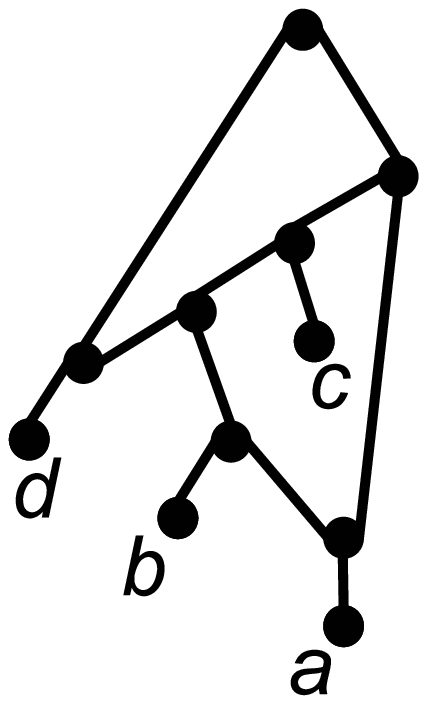}
& ~~~~~~
&\includegraphics[scale=0.3]{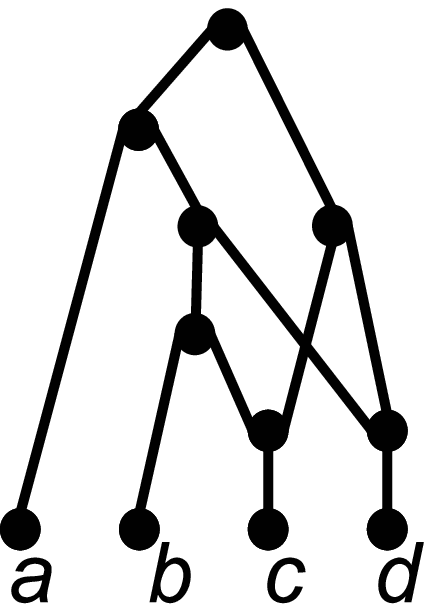} 
&\includegraphics[scale=0.3]{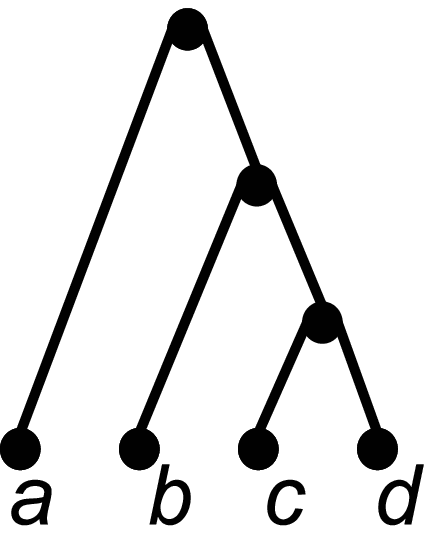}\\
(a) & & (b) &(c)
\end{tabular}
\caption{(a) A level-2 network $N$ for which $\mathcal S(N)$
is not a weak hierarchy. The phylogenetic network $N$ depicted in (b) does not display the
phylogenetic tree $T$ depicted in (c) but
$\mathcal C(T)$ is tree-consistent with $N$. }
\label{FigLevel2WeakHierarchy}
\end{figure}

We conclude this section with remarking that as the example of 
the level-$2$ network $N$ presented in 
Fig.~\ref{FigLevel2WeakHierarchy}(a)  shows, the result analogous 
to Proposition~\ref{softwired-weak-hierarchy} does not
hold for level-$2$ networks since 
$\{\{a,b,c\},\{a,b,d\},$
$\{b,c,d\}\}\subseteq\mathcal S(N)$ but 
$\{a,b,c\} \cap \{a,b,d\} \cap \{b,c,d\}=\{b\}$.
Furthermore, we remark that  as the example
of the level-2 network $N$ depicted in
 Fig.~\ref{FigLevel2WeakHierarchy}(b) combined with
the cluster system $\mathcal C(T)$ induced by the phylogenetic tree 
$T$ depicted in Fig.~\ref{FigLevel2WeakHierarchy}(c) shows, 
a cluster system $\mathcal C(T)$ induced by a phylogenetic tree $T$
can be contained in the cluster system $\mathcal S(N)$ of a level-$2$
network $N$ and $N$ need not display $T$.

\section{Simple level-$1$ Networks}\label{Simple-level-$1$}

%Then we call for some
%arc $uv\in A(G)$, the vertex $u$ a {\em father} of $v$ and 
%$v$ a {\em son} of $u$.
In this section we turn our attention 
to studying simple level-$1$ networks.
In particular, we establish a fundamental property of these
networks with regards to encodings of level-$1$ networks. 
To do this, we require some more definitions. We start with 
the definition of the triplet system $\mathcal R(N)$ induced
by a phylogenetic network $N$.
 
\begin{figure}[!hbt]
\begin{center}
\begin{tabular}{ccc}
\includegraphics[scale=0.45]{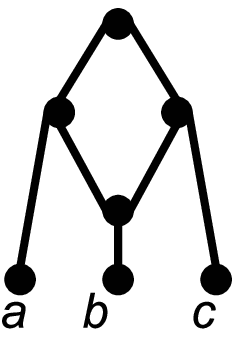} &
\includegraphics[scale=0.45]{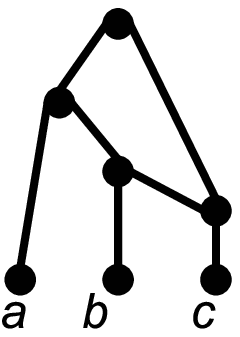} &
\includegraphics[scale=0.45]{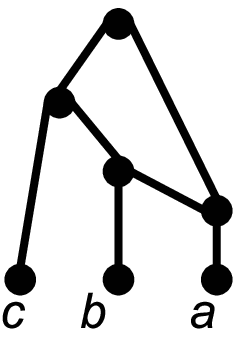} \\
$SL_1^1(T)$ & $SL_1^2(T)$ & $SL_1^3(T)$
\end{tabular}
\end{center}
\caption{\label{simple-level-$1$} The three non-isomorphic 
simple level-$1$ networks 
on $\{a,b,c\}$ that all reflect the triplet 
system $\mathcal R=\{a|bc, c|ab\}$. 
}
\end{figure}

Suppose $N$ is a phylogenetic network. If $Y\subseteq X$
is a subset of $X$ of size 3, then $N$ induces a triplet $t$ on $X$
 by taking
$t$ to be a minimal subtree of $N$ with leaf set $Y$
and suppressing resulting degree two vertices of $t$. 
The set of triplets induced on $X$ by $N$ this
way is the triplet system $\cR(N)$. Two properties of this 
triplet systems should be noted. First, every triplet in $\cR(N)$
is consistent with $N$, where a triplet $x|yz$ is called 
{\em consistent} with a phylogenetic network $N$ if $x,y,z\in X$
 and there exist two vertices $u,v\in V(N)$  and pairwise internally
vertex-disjoint paths in $N$ from $u$ to $y$, $u$ to $z$, $v$ to 
$u$ and $v$ to $x$.
Note that a triplet system $\cR$ is called {\em consistent} 
with a phylogenetic network $N$ if every triplet in $\mathcal R$
is consistent with $N$. For convenience, we will sometimes
say that a phylogenetic network $N$ is consistent with a
triplet $t$ (or a triplet system $\cR$) if $t$ (or $\cR$) is  consistent
with $N$. In case $\cR$ is consistent with a phylogenetic
network $N$ and $\mathcal R=\cR(N)$ then  we say that 
$\mathcal R$ {\em reflects} $N$. Alternatively, we will say that
$\mathcal R$ is {\em reflected} by $N$. For example,
the triplet set $\mathcal R=\{a|bc, c|ab\}$
is reflected by the three simple level-$1$ networks 
$SL_1^i(T)$, $i\in\{1,2,3\}$ on  $\{a,b,c\}$ 
depicted in Fig.~\ref{simple-level-$1$} which appeared in
slightly different form in \cite{JNS2006}.

Second, the triplet system $\cR(N)$ is always dense, 
where a triplet system $\cR$ on $X$ is called {\em dense} if
for any three elements in $a,b,c\in X$ there exits a triplet $t\in \cR$
such that $L(t)=\{a,b,c\}$. Arguably unassumingly looking, the
concept of a dense triplet set has proven vital for level-$k$
network reconstruction, $k\geq 1$, from triplet systems. More precisely,
 the only known polynomial time algorithms 
for constructing level-$1$ and level-2 networks $N$  
consistent with such triplet systems construct $N$, (if it exist)
by essentially building it up recursively from 
simple level-$1$ and simple level-$2$ networks 
\cite{JanssonSung2006,IKKS2008}. If the assumption
that $\mathcal R$ is dense is dropped however, then it is NP-hard
to decide if there exists a level-$k$ network, $k=1,2$, consistent with 
$\mathcal R$ \cite{JanssonSung2006,IKKS2008}. 
For larger values of $k$, a polynomial time algorithm 
for constructing a level-$k$ network from a dense triple set
was recently presented in \cite{ToHabib2009}.

The next result is rather technical\footnote{ 
A case analysis based alternative proof of this 
result may be found in~\cite{GBP2008}.}
but plays a crucial role in the proof of 
our main result (Theorem~\ref{theorem-montpellier}) as it shows that
although all three simple level-$1$ networks depicted
in Fig.~\ref{simple-level-$1$} reflect the same triplet set this property is
lost when adding an additional leaf to a non-cut-arc of each of them. 
For a directed graph $G$ these arcs are the elements in $A(G)$
whose removal disconnect $G$.  To establish our result,
we require some more definitions and
notations.

Suppose $N$ is a phylogenetic network and $a,b\in V(N)$
such that $a$ is below $b$. If $c$ is a further vertex in $V(N)$ and
$a \preceq_N b$ and $c\preceq_N b$ holds
then we call $b$ a {\em common ancestor} of $a$ and $c$.
A {\em lowest common ancestor} $lca_N(a,c)$ of $a$ and $c$
is a common ancestor of $a$
and $c$ and no other vertex below $lca_N(a,c)$ is a common ancestor 
of $a$ and $c$. 
Note that in a level-0 or level-$1$ network $N$, the lowest common 
ancestor between any two distinct leaves of $N$ is always unique
whereas this need not be the case for level-$k$ networks with larger $k$.

Now suppose $N$ is
one of the simple level-$1$ networks $SL_1^i(T)$, $i\in \{1,2,3\}$,
on $X=\{a,b,c\}$ depicted in Fig.~\ref{simple-level-$1$}. Let $e=uv\in  A(N) $
be a non-cut arc and suppose that $d\not\in X$. Then we denote
by $N_e\oplus d$ the level-$1$ network obtained from $N$ 
by adding a new vertex $w$ to $V(N)$ and replacing
$e$ by the arcs $uw$, $wv$, and $wd$. We remark that if the knowledge of $e$ is of no 
relevance to the presented argument, then we
will write $N\oplus d$ rather than  $N_e\oplus d$.

\begin{lemma}\label{lemma-montpellier}
Suppose $X=\{a,b,c,d\}$ and $T=\{a|bc, c|ab\}$. Then, 
for any two distinct $i,j\in \{1,2,3\}$,
$$
\cR(SL_1^i(T)\oplus d)\not=\cR(SL_1^j(T)\oplus d).
$$   
\end{lemma}
\pf Put $N^k:=SL_1^k(T)$, $k\in\{1,2,3\}$, and assume that 
there exist distinct $i,j\in\{1,2,3\}$ 
and non-cut-arcs $e_i\in A(N_i)$
and $e_j\in A(N_j)$ such that 
$\cR(N^i_{e_i}\oplus d)=\cR(N^j_{e_j}\oplus d)$. 
By symmetry, it suffices to consider
the cases $(i,j) \in \{(2,1),(2,3)\}$.
For $k \in \{1,2,3\}$, let $u_k,v_k\in V(N^k)$ such that $e_k=u_kv_k$.
Also for $k \in \{1,2,3\}$, let
 $w_k\not\in V(N^k)$ denote the new vertex in $V(N^k_{e_k}\oplus d)$
such that
by replacing the arc $e_k$ by the
arcs $u_kw_k$, $w_kv_k$, and adding the arc $w_kd$ the new
network $N^k_{e_k}\oplus d$ is obtained from $N^k$.
Note that for all $k\in\{1,2,3\}$, both $N^k$ and $N^{k}\oplus d$
have the same root and the same
reticulation vertex which we denote by $\rho_k$ and $r_k$, respectively.
Furthermore, for all $x,y \in \{a,b,c\}$ we have
$lca_{N^k}(x,y) = lca_{N^k\oplus d}(x,y)$
We distinguish the cases that $u_2=\rho_2$ and 
that $u_2\not=\rho_2$.

Suppose first that $u_2=\rho_2$ and put $l=lca_{N^2}(a,b)$. 
Then $e_2\in \{\rho_2 r_2, \rho_2 l\}$.
We first establish that $j\not=1$.
Assume for contradiction that $j=1$. 
%Then $\cR(N^2_{e_2}\oplus d)=\cR(N^1_{e_1}\oplus d)$.
For all $s,t\in V(N^2_{e_2}\oplus d)$ such that $t$ is below $s$ 
 denote a path from $s$ to $t$ in $N^2_{e_2}\oplus d$ by $P_{st}$.
Observe that $d|ac\in \cR(N^2_{e_2}\oplus d)$, 
$e_2\in \{\rho_2r_2, \rho_2 l\}$, holds. Indeed,
since $c|ab\in \cR(N^2)$ the paths $P_{la}$ and 
$P_{lc}$ exist and do not have an internal vertex in common.
Furthermore, $w_2\not=l$ and either the arc $w_2l$
or the arcs $\rho_2 l$ and $\rho_2 w_2$
exist. In both cases the paths $P_{w_2l}$, $P_{\rho_2 l}$ and $P_{\rho_2 w_2}$,
consisting of the arcs $w_2l$, $\rho_2 l$ and $\rho_2 w_2$, respectively,
do not have an internal vertex in common with either
$P_{la}$ or $P_{lc}$. Thus, $d|ac\in \cR(N^2_{e_2}\oplus d)$, as required.
By assumption,  $d|ac\in\cR(N^1_{e_1}\oplus d)$ follows which is impossible
since $lca_{N^1}(a,c)=\rho_1$ and so $d|ac=\cR(N^1_e\oplus d)$,
for all non-cut arcs $e\in A(N^1)$.
Thus, $j\not=1$, as required.

If $j=3$ then $\cR(N^2_{e_2}\oplus d)=\cR(N^3_{e_3}\oplus d)$ and
$v_2=r_2$ or  $v_2=l$. If $v_2=r_2$ then
$b|cd\in \cR(N^2_{e_2}\oplus d)=\cR(N^3_{e_3}\oplus d)$ follows.
But then $w_3$ cannot be a vertex on the path in $N^3_{e_3}\oplus d$
from $\rho_3$ to $b$ or on  the path from $\rho_3$ to $r_3$ 
which avoids $lca_{N^3}(a,c)$
Thus, $u_3=lca_{N^3}(a,b)$ and so 
$b|da\in \cR(N^3_{e_3}\oplus d)=\cR(N^2_{e_2}\oplus d)$ which is impossible
as $lca_{N^2\oplus d}(a,d)= \rho_2$ and thus always above 
$l$. If
$v_2=l$ then $d|bc,c|bd\in \cR(N^2_{e_2}\oplus d) 
=\cR(N^3_{e_3}\oplus d)$ 
follows which is again impossible since if $w_3$ lies on the path
from $lca_{N^3}(a,c)$ to $r_3$ then
$d|bc\not \in \cR(N^3_{e_3}\oplus d)$ and if not then 
$u_3=\rho_3$
and so $c|bd\not\in \cR(N^3_{e_3}\oplus d)$.  Thus, $j\not=3$.

Now suppose that $u_2\not=\rho_2$. Then $u_2\in \{l,lca_{N^2}(b,c)\}$
Observe that arguments similar 
to the previous ones imply that
$a|cd,a|bd,c|ad,c|bd\in \cR(N^2_{e_2}\oplus d)$ holds for all
$u\in \{l,lca_{N^2}(b,c) \}$. If $j=1$ then 
$u_1\not=\rho_1$ as otherwise  $a|bd$ or $c|bd$ 
does not belong to $\cR(N^1_{e_1}\oplus d)=\cR(N^2_{e_2}\oplus d)$.
Thus $v_1=r_1$ and 
$u_1\in\{lca_{N^1}(b,a),lca_{N^1}(b,c)\}$. If 
$u_1=lca_{N^1}(b,a)$ then $a|cd\not\in  \cR(N^1_{e_1}\oplus d)$
which is impossible. Swapping the roles of $a$ and $c$ in the 
previous argument shows that $u_1=lca_{N^1}(b,c)$
cannot hold either. Thus, $j\not=1$. 

If $j=3$ then again since $a|cd,c|ad\in \cR(N^1_{e_1}\oplus d)$,
it follows  that $e_3$
must be an arc on the path $P$
from $lca_{N^3}(a,c)$ to 
$r_3$. Note that similar arguments as the ones 
used above imply that either $d|bc$ or $b|cd$ is contained in
$\cR(N^2_{e_2}\oplus d)$. But 
$b|cd\not\in \cR(N^3_{e_3}\oplus d)=\cR(N^2_{e_2}\oplus d)$
and so $d|bc\in \cR(N^2_{e_2}\oplus d)$ must hold.
But this is impossible since then $c|ad,d|bc\in \cR(N^2_{e_2}\oplus d)$
but there exists no non-cut-arc $e$ on $P$ 
such that both triplets are simultaneously contained in 
$\cR(N^3_e\oplus d)$. Thus, $j\not =3$.
\epf

\section{Encodings of Level-$k$ Networks}\label{alternativeEncoding}
\label{tripletEncoding}
\label{consistency-level-$1$}

In this section, we characterize those level-$1$ networks $N$ 
that are encoded by the triplet system $\mathcal R(N)$, or equivalently
the tree system $\mathcal T(N)$, or equivalently the 
cluster system $\mathcal S(N)$ they induce. In addition, we present an example that
illustrates that our arguments cannot be extended to establish the
corresponding result for level-2 networks and therefore to
level-$k$ networks with $k\geq 3$.

Bearing in mind that there exist triplet system  
which can be reflected by more than 
one level-$1$ network, we denote the collection of all level-$1$ networks 
that reflect a triplet system $\mathcal R$ by
$\mathfrak{L}_1(\mathcal R)$. Clearly, if $\mathcal R$
is reflected by a level-$1$ network $N$ then $N\in \mathfrak{L}_1(\cR(N))$
and so $|\mathfrak{L}_1(\cR(N))|\geq 1$. Similarly, we 
denote for a tree system $\mathcal T$ the collection
of all level-$1$ networks $N$ for which  $\mathcal T=\mathcal T(N)$ holds by
$\mathfrak{L}_1(\mathcal T)$, and for a cluster system $\mathcal C$ the
collection of all level-$1$ networks $N$ for which  
$\mathcal C=\mathcal S(N)$ holds by $\mathfrak{L}_1(\mathcal C)$. 
As in the case of triplet systems, 
there exist tree systems $\mathcal T$ and cluster systems $\mathcal C$ 
with $|\mathfrak{L}_1(\mathcal T)|\geq 1$ and 
$|\mathfrak{L}_1(\mathcal C)|\geq 1$, respectively. 

Clearly, any cluster $C\subseteq X$ induces a triplet system
$\mathcal R(C)$ of triplets on $X$ defined by putting 
$$
\mathcal R(C)=\{c_1c_2|x : c_1,c_2\in C \mbox{ and } x\in X-C\}.
$$
Thus, any non-empty cluster system $\mathcal C$ on $X$
induces a triplet system $\mathcal R(\mathcal C)$ defined by 
putting
$\mathcal R(\mathcal C):=\bigcup_{C\in\mathcal C} \mathcal R(C)$.
The next result establishes a link between the triplet system induced
by a level-$1$ network $N$ and the triplet system 
$\mathcal R(\mathcal S(N))$.

\begin{lemma}\label{new}
Suppose $N$ is a level-$1$ network with at least 3 leaves.
Then 
$$
\mathcal R(N)=\bigcup_{T\in \mathcal T(N)}\mathcal R(T)=
\bigcup_{C\in \mathcal S(N)}\mathcal R(C).
$$
\end{lemma}
\pf That $\bigcup_{T\in \mathcal T(N)}\mathcal R(T)=
\bigcup_{C\in \mathcal S(N)}\mathcal R(C)$ holds is trivial.
Also it is straight forward to see that 
$\bigcup_{T\in \mathcal T(N)}\mathcal R(T)\subseteq\mathcal R(N)$.
To see the converse set inclusion, suppose that $t\in \mathcal R(N)$.
Let $x_1,x_2,x_3\in X$ such that 
$t=x_1x_2|x_3$. Then with $lca(x_1,x_2):=lca_N(x_1,x_2)$
we have $x_3\not\in C_N( lca_N(x_1,x_2))$ and
$lca(x_1,x_2)$ does not equal the root $\rho_N$ of 
$N$. Let $P_i$ denote a path from 
$\rho_N$  to $x_i$, $i=1,2$ and let $T$ denote
the phylogenetic tree on $X$ obtained from $N$ by modifying
all reticulation vertices $v$ of $N$ in the following way.
If $v\not\in V(P_1)\cup V(P_2)$
then randomly delete one of the incoming arcs of $v$ and suppress the 
resulting degree 2 vertex. If this results in the decrease of the outdegree
of the root $\rho_N$ of $N$  then identify $\rho_N$ with is unique child.
If $v\in V(P_i)$, $i=1,2$, then delete that incoming arc of 
$v$ that is not an arc of $P_i$ and suppress the resulting degree
2 vertex. Clearly, $T$ is displayed by $N$ and so 
$t\in \mathcal \bigcup_{T\in \mathcal T(N)}\mathcal R(T)$.
Thus, $\mathcal R(N)\subseteq \bigcup_{T\in \mathcal T(N)}\mathcal R(T)$
must hold which implies the lemma.
\epf

\begin{figure}[!hgt]
\centering
\includegraphics[scale=0.3]{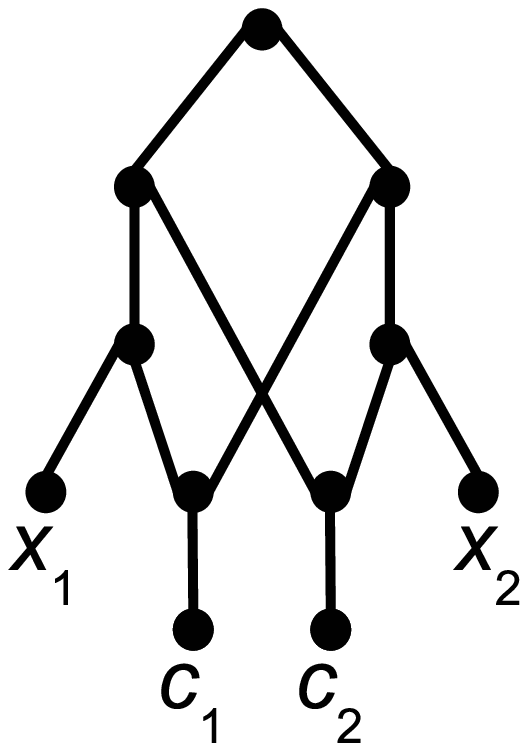}
\caption{A level-2 phylogenetic network $N$ with 
$c_1c_2|x_1\in \mathcal R(N)$, but $\{c_1,c_2\} \not\in \mathcal S(N)$.
}
\label{FigLevel2Counter-changed}
\end{figure}

Note that as the example of the level-$2$ network depicted in
Fig.~\ref{FigLevel2Counter-changed} shows,
%Fig.~\ref{FigLevel2Counter-changed} shows,
the relationship between the triplet system of a level-$1$ network
$N$ and the triplet system induced by the clusters in $\mathcal S(N)$
%exploited in the proof of Lemma~\ref{new} 
does not hold for level-$2$
networks. 
%
%\begin{figure}[!hgt]
%\centering
%\includegraphics[scale=0.3]{Level2Counterexample.eps}
%\caption{A level-2 phylogenetic network $N$ with 
%$c_1c_2|x_1\in \mathcal R(N)$, but $\{c_1,c_2\} \not\in \mathcal S(N)$.
%}
%\label{FigLevel2Counter-changed}
%\end{figure}

To prove the main result of this note (Theorem~\ref{theorem-montpellier})
which we will do next, we require some additional 
definitions and notations. Suppose $N$ is a phylogenetic network.
Then we call a subset $\{x,y\}\subseteq  X$ a 
{\em cherry} of $N$ if there exists a vertex $v\in V(N)$ 
such that $vx,vy\in A(N)$. Furthermore, if $N$ is a level-$1$
network and $x\in X$ then
we denote by $N-x$ the level-$1$ network obtained from $N$ 
by removing $x$ (and its incident arc) and suppressing 
the resulting degree 2 vertex. In addition, we say that
$N$ is a {\em strict} level-$1$ network if $N$ is not a
phylogenetic tree. Finally, to a triplet system $\mathcal R$
and some $x\in \bigcup_{t\in \mathcal R}L(t)$, we associate the
triplet set $\mathcal R_x:=\{t\in \mathcal R : x\not\in L(t)\}$.

Armed with these definitions and notations we are now ready to 
establish our main result.

\begin{theorem}\label{theorem-montpellier}
Suppose $N$ is a  level-$1$ network with at least 3 leaves. 
Then the following statements are equivalent
\begin{enumerate}
\item[(i)] $N$ contains a blob with four vertices.
\item[(ii)]  $|\frak{L}_1(\mathcal R (N))|>1$.
\item[(iii)]  $|\frak{L}_1(\mathcal S (N))|>1$.
\item[(iv)]  $|\frak{L}_1(\mathcal T (N))|>1$.
\end{enumerate}
\end{theorem}
\pf
(i) $\Rightarrow$ (iv):  This is an immediate consequence
of the fact that all simple level-$1$ networks depicted
in Fig.~\ref{simple-level-$1$} induce the same set of phylogenetic trees.

\noindent (iv) $\Rightarrow$ (ii): 
Suppose that $N$ is a level-$1$ network such that 
$|\frak{L}_1(\mathcal T (N))|>1$. Then there exists a level-$1$ network
$N'$ distinct from $N$ such that $\mathcal T(N)=\mathcal T (N')$. 
Combined with Lemma~\ref{new}, 
$
\mathcal R(N)=\bigcup_{T\in \mathcal T(N)} \mathcal R(T)=
\bigcup_{T\in \mathcal T(N')} \mathcal R(T)=\mathcal R(N')
$
follows and so  $N' \in \frak{L}_1(\mathcal R (N))$. 
Thus,  $|\frak{L}_1(\mathcal R (N))|>1$, as required.

\noindent (ii) $\Rightarrow$ (i)
We will show by induction on the number $n$ of leaves of $N$
that if every blob in $N$ contains at least 5 vertices 
then $|\frak{L}_1(\cR(N))|=1$. Suppose $N$ is a level-$1$ network
with $n$ leaves such that every blob of $N$ contains at least 5 vertices.
Note that we may assume that $N$ contains at least one  blob since otherwise
$N$ is a phylogenetic tree and so $|\frak{L}_1(\cR(N))|=1$
clearly holds. 
But then $n\geq  4$. If $n=4$ then, using Lemma~\ref{lemma-montpellier},
it is straightforward to verify that $|\frak{L}_1(\cR(N))|=1$. 
\old{
Suppose first that $N$ has four leaves 
$a,b,c,d$, say, and assume without loss of generality
that $T'=\{a|bc, c|ab\} \subseteq \cR(N)$.
Then there exist
some $i\in \{1,2,3\}$ such that $N$ equals $SL_1^i(T' )\oplus d$.
But then, by Lemma~\ref{lemma-montpellier},  $|\frak{L}_1(\cR(N))|=1$ 
follows.
}

Suppose $n> 4$. Assume for every level-$1$ network $N_0$ 
with $n_0<n$ leaves that $|\frak{L}_1(\cR(N_0))|=1$ holds
whenever $N_0$ is a phylogenetic tree or 
every blob in $N_0$ contains at least 5 vertices. Suppose for contradiction
that $|\frak{L}_1(\cR(N))| \geq 2$. Choose some  $N'\in \frak{L}_1(\cR(N))$
distinct from $N$.
Then $\mathcal R:=\cR(N)=\cR(N')$. We distinguish
the cases that $N$ contains a cherry and that it does not.

Suppose first that $N$ contains a cherry $\{x,y\}$. 
Without loss of generality, we may assume that this cherry 
is as far away from the root of $N$ as possible.
Then since $N$
is a strict level-$1$ network all of whose blobs contain
at least 5 vertices, $N-x$ must enjoy the same property with 
regards to its blobs (if $N-x$ still has blobs). But then, by induction 
hypothesis, $|\frak{L}_1(\cR(N-x))|=1$ and so $N-x$
is the unique level-$1$ network that reflects  
$\cR(N-x)=\mathcal R_x$. Since by the choice of $x$, for every leaf $z$ in $N$ 
distinct from $x$ and $y$, only the triplet $z|xy$ out of the 3 possible
triplets on $\{x,y,z\}$ is contained in $\mathcal R=\cR(N')$, it follows
that $\{x,y\}$ must also be a cherry in $N'$. But then 
$N=N'$ which is impossible.
Thus, $|\frak{L}_1(\mathcal R(N))|=1$ must hold in this case.

Now suppose that $N$ does not contain a cherry. Then there exists 
a blob $B$ in $N$ such that all cut-arcs that start with a vertex in $B$
must end in a leaf of $N$. For each such leaf $z$, which we will also
call a leaf of $B$, we denote by $z'$
the vertex of $B$ such that $z'z$ is that cut-arc of $N$. Furthermore,
denote by $p$ the leaf of $B$ such that $p'$ is the reticulation vertex in
$B$. Let $y_1$ and $y_2$ the vertices in $V(N)-V(B)$ such that $y_1'$ and $y_2'$
are  the two parent vertices of 
$p'$ in $B$. Note that the root $\rho=\rho_B$ of $B$ 
could be $y_1'$ or $y_2'$ but not both and that whenever
$y_i'\not=\rho$, $i=1,2$, then $y_i$ is a leaf of $B$ (hence
the abuse of notation). Without 
loss of generality, we may assume that the path $P_{\rho {y_1'}}$
from $\rho$ to $y_1'$ in $B$ is at least as long as the path $P_{\rho {y_2'}}$
 from $\rho$ to $y_2'$ in $B$ (where we allow paths 
of length zero). Thus, $y_1$ must be a leaf of $B$.
Since $P_{\rho {y_1'}}$ is at least
as long as $P_{\rho {y_2'}}$ and, by assumption on $N$, $B$
contains at least 5 vertices,
there must exist a leaf $y$ of $B$
distinct from $y_1$  such that $y'\in V(P_{\rho {y_1'}})$. 
Note that we may assume without loss of generality
that $y'$ is the predecessor of $y_1'$ on that path.
 We distinguish the cases that $|V(B)|>5$ and that $|V(B)|=5$.

Suppose first $|V(B)|>5$.
Since a blob in the level-$1$ network
$N-y_1$ has clearly at least $5$ vertices, we have
$|\frak{L}_1(\cR(N-y_1))|=1$  by the induction hypothesis. But then 
 $N-y_1$ is the unique level-$1$ network that reflects 
$\cR(N-y_1)=\mathcal R_{y_1}$.
 Consequently, since $ \mathcal R=\cR(N')$ we have  $N'-y_1=N-y_1$. 
To see that $N$ equals $N'$ suppose $z$ is a leaf of 
$B$ distinct from $y_1,y,p$ (which must exist by assumption on $B$). Then
either $t:=z|yp,p|yz\in \mathcal R$ or $t,y|zp\in \mathcal R$ holds.
We only discuss the case that $t,p|yz\in \mathcal R$ since the case 
$t,y|zp\in \mathcal R$
is symmetric. 
%Since $y_1$ is neither a leaf of $z|yp$ nor of $p|yz$, 
%both triplets must be contained in $\cR(N')$. 
Let $B^-$ denote the
blob in $N-y_1$ obtained from $B$
by deleting $y_1$ plus its incident arc and suppressing the resulting degree
2 vertex. Since $z$, $y$, and $p$ are leaves of 
$B^-$ and the choice of $y_1$ implies that 
$y_1y|p, y|y_1p\in \cR(N)=\cR(N')$, it follows that there exists some
blob $B'$ in $N'$ such that 
$B^-=B'-y_1$.  Moreover, the suppressed degree 2 vertex of
$V(B')$ is adjacent (in $B'$) with $y'$ and 
$p'$, respectively, since
otherwise $y_1|yp\in \cR(N')=\cR(N)$ would hold which contradicts
the choice of $y_1$. Thus  $N=N'$ and so $|\frak{L}_1(\cR(N))|=1$
must hold in case $|V(B)|>5$.

We conclude with analyzing  the case $|V(B)|=5$. 
Then either $\rho=y_2$ and so $B$ has, in addition to  
the leaves $y_1,y,p$, precisely one more leaf $z$, or
$\rho\not=y_2$ and the leaves of $B$ are $y_1,y_2,y$ and $p$.
We first consider the case $\rho\not=y_2$. 
Consider the level-$1$ network $N-\{y_1,y_1'\}$ obtained from $N$ 
by removing $y_1$, its parent vertex $y_1'$ and their 3 incident
arcs (plus suppressing resulting degree 2 vertices)
thus effectively turning $B$ into a phylogenetic tree on
the leaves $y,p,y_2$, i.e. the
triplet $t:= y|py_2$. Put $\mathcal R^{t}:=\mathcal R_{y_1}\cup\{t\}$.
Since either $N-\{y_1,y_1'\}$ is a phylogenetic tree or a strict level-$1$
network such that each of its blobs contains at least 5 vertices, the
induction hypothesis implies 
$|\frak{L}_1(\cR(N-\{y_1,y_1'\}))|=1$. Thus,
$N-\{y_1,y_1'\}$ is the unique level-$1$ network that 
reflects $\mathcal R^{t}$. 
Note that the only way to turn $N-\{y_1,y_1'\}$ into a level-$1$ network that,
in addition to reflecting $\mathcal R^{t}$, is also consistent with $t':=y_2|py\in\mathcal R$
is to replace $t$ by one of the level-$1$ networks $SL_1^j(\{t,t'\})$,
$j\in Y:=\{1,2,3\}$. Denote that element in $Y$ by  $j_{N}$. Since 
$\cR(N)=\cR(N')$  it follows that the level-$1$ network  obtained
from $N'$ by removing $y_1$, its parent vertex, and their 3 incident
arcs (suppressing resulting degree 2 vertices)
must equal $N-\{y_1,y_1'\}$ with $t$ replaced by one of  $SL_1^j(\{t,t'\})$,
$j\in Y$. Denote that element in $Y$ by  $j_{N'}$.
 Since $\{ y_1|py_2, y_2|py_1, y_2|y_1y, p|y_1y, y|y_1p, y_2|py, t\}
\subseteq \mathcal R=\cR(N')$ it is easy to check 
that $j_N=j_{N'}$ must hold and so
$N$ and $N'$ must be 
equal which is again impossible. 
Thus, $|\frak{L}_1(\cR(N))|=1$ must hold in case 
$\rho\not=y_2$. Using arguments similar to the previous ones
it is straight-forward $N=N'$ and thus $|\frak{L}_1(\cR(N))|=1$
must hold in case $\rho=y_2$. 

\noindent (iv) $\Rightarrow$ (iii): Suppose that $N$ is a level-$1$ network
with $\frak{L}_1(\mathcal T (N))|>1$. Then
there exists a level-$1$ network
$N'\in \frak{L}_1(\mathcal T (N))$ distinct from $N$ with
$\mathcal T(N)=\mathcal T(N')$. But then
$\mathcal S(N)=\bigcup_{T \in \mathcal T (N)} \mathcal C(T)=
\bigcup_{T \in \mathcal T (N')} \mathcal C(T)=
\mathcal S (N')$ and so $N' \in \frak{L}_1(\mathcal S (N))$. Thus,
$\frak{L}_1(\mathcal S (N))|>1$.

\noindent (iii) $\Rightarrow$ (ii):
Suppose that $N$ is a level-$1$ network with $|\frak{L}_1(\mathcal S (N))|>1$.
Then there exists a level-$1$ network
$N' \in \frak{L}_1(\mathcal S (N))|>1$ distinct from $N$ such that 
$\mathcal S (N)=\mathcal S (N')$. But then Lemma~\ref{new}
implies 
$
\mathcal R(N)=\bigcup_{C\in \mathcal S(N)} \mathcal R(C)=
\bigcup_{C\in \mathcal S(N')} \mathcal R(C)=\mathcal R(N')
$
and so $N'\in \frak{L}_1(\mathcal R (N))$. Hence, 
$\frak{L}_1(\mathcal R(N))|>1$.
\epf

It should be noted that Theorem~\ref{theorem-montpellier} immediately implies

\begin{corollary}
Let $N$ be a level-$1$ network with at least 3 leaves. The number of 
non-isomorphic level-$1$ networks $N'$
that reflect $\mathcal R(N)$ (or equivalently for which
$\mathcal T(N)=\mathcal T(N')$ or equivalently
$\mathcal S(N)=\mathcal S(N')$ holds) is $3^b$, where
$b$ is the number of blobs of $N$ of size four.
\end{corollary}

\begin{figure}[!hgt]
\centering
\includegraphics[scale=0.3]{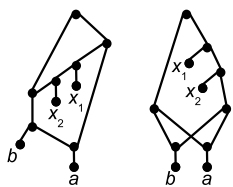}
\caption{Both simple level-$2$ networks reflect the triplet set
$\{a|x_1b,$ $b|x_1a,$ $x_1|ab,$ $a|x_2b,$ $b|x_2a,$ $x_2|ab,$ $x_1|x_2a,$ 
$a|x_1x_2,$ $x_1|x_2b,$ $b|x_1x_2\}$.
}
\label{FigLevel2}
\end{figure}

We remark that the strategy 
underlying the proof of Theorem~\ref{theorem-montpellier}
does not immediately extend to level-$k$ networks with $k\geq 2$. The main
reasons for this are that, as already mentioned above, for $k\geq 2$ the 
number of distinct level-$k$ generators grows 
exponentially in $k$ \cite{GBP2009}.
Also the problem of understanding when two distinct 
simple level-$2$ networks reflect the same set of
triplets is far less well understood. For example, 
consider the two level-$2$ networks
depicted in Figure~\ref{FigLevel2}. Each one of them is a simple level-$2$
network obtained by hanging leaves of the sides of 
the level-$2$ generators $\mathcal G^2_a$ and $\mathcal G^2_b$ 
depicted in Figure~\ref{level-k-generators}. 
As can be quickly
verified, both networks reflect the same triplet set. However adding
additional leaves to both networks by subdividing
the arc one of whose end vertices forms an arc with $x_1$ and the
other forms an arc with $x_2$ and attaching additional leaves 
results in two distinct level-$2$ networks that still reflect  
the same triplet system.  Regarding the accurate reconstruction of
level-$k$ networks from e.g. triplet data, this results highlights  
a serious limitation of level-$2$ networks (and probably
level-$k$ networks in general) as two such network with very different
structure might reflect the same triplet set.

We conclude with remarking that 
phylogenetic trees on $X$  can also be viewed as
trees together with a bijective labelling map between $X$ and
the leaf set of such trees. Taking this point of view,
phylogenetic trees were generalized in
\cite{HuberMoulton2006}  to {\em MUL-trees} 
by allowing two or more leaves of that tree to have the
same label. For example, the tree obtained
from the phylogenetic tree depicted
in Figure~\ref{FigLevel1}(c) by replacing the leaf labelled $a$ 
by the cherry labelled $\{a,b\}$ is such  a tree. In fact,
this is the MUL-tree induced by
the level-$1$ network  $N$ depicted in Figure~\ref{FigLevel1}(a)
that shows all paths from the root  of $N$ to all leaves of $N$.
For a level-$1$ network $N$ it is easily seen that the MUL-tree $\mathcal M(N)$
induced by $N$ this way is in fact an encoding of $N$ in the sense that 
$N$ is the unique level-$1$ 
network that can give rise to $\mathcal M(N)$.

%\begin{acknowledgements}
%\end{acknowledgements}

\bibliographystyle{alpha}
\bibliography{Phylogeny}   % name your BibTeX data base

\end{document}